\documentclass[10pt,amssymb]{article}
\usepackage{latexsym,amsfonts,amsmath}

\title{Riemann minimal surfaces in higher dimensions}

\author{S. Kaabachi and F. Pacard}

\begin{document}
\maketitle
% LES NOUVELLES COMMANDES:

       \newcommand{\RS}{\mathbb{R} \times S^{n-1}}
       \newcommand{\Sn}{S^{n-1}}

       \newcommand{\calL}{\mathcal L}
       \newcommand{\calH}{\mathcal H}
       \newcommand{\calQ}{\mathcal Q}

       \newcommand{\tet}{y}
       \newcommand{\e}{\varepsilon}

       \newcommand{\vph}{\varphi}

       \newtheorem{theorem}{Theorem}
       \newtheorem{proposition}{Proposition}
       \newtheorem{corollary}{Corollary}
       \newtheorem{lemma}{Lemma}
       \newtheorem{definition}{Definition}
       \newtheorem{remark}{Remark}

%LES NOUVEAUX ENVIRONNEMENTS :
       \newenvironment{dem} {\noindent {\bf Proof\/}: } {\null \hfill
$\Box$}
       \newenvironment{case} {\noindent {\textit{Case}\/} }{\null}
       \newenvironment{step} {\noindent {\textit{Step}\/} }{\null}
       \renewcommand{\theequation}{arabic\thesection.\arabic{equation}}
       \def\thesection{\arabic{section}}
       \def\thetheorem{\arabic{section}.\arabic{theorem}}
       \def\theproposition{\arabic{section}.\arabic{proposition}}
       \def\thecorollary{\arabic{section}.\arabic{corollary}}
       \def\thelemma{\arabic{section}.\arabic{lemma}}
       \def\thedefinition{\arabic{section}.\arabic{definition}}
       \def\theremark{\arabic{section}.\arabic{remark}}

%LES NUMEROTATIONS :
       \numberwithin{equation}{section}
       \numberwithin{proposition}{section}
       \numberwithin{corollary}{section}
       \numberwithin{lemma}{section}
       \numberwithin{remark}{section}
       \numberwithin{definition}{section}

{\bf Abstract.} We prove the existence of a one parameter family of
minimal embedded hypersurfaces in $R^{n+1}$, for $n \geq 3$, which
generalize the well known $2$ dimensional "Riemann minimal
surfaces". The hypersurfaces we obtain are complete, embedded,
simply periodic hypersurfaces which have infinitely many parallel
hyperplanar ends. By opposition with the $2$-dimensional case, they
are not foliated by spheres.

\medskip

{\bf R\'esum\'e.} Nous prouvons l'existence d'une famille \`a un
param\`etre d'hypersurfaces  de $R^{n+1}$, pour $n\geq 3$, qui sont
minimales et qui g\'en\'eralisent les surfaces minimales de Riemann.
Les hypersurfaces que nous obtenons sont des hypersurfaces
compl\`etes, simplement p\'eriodiques et qui ont une infinit\'e de
bouts hyperplans parall\`eles. Contrairement au cas des surfaces,
i.e. $n=2$, ces hypersurfaces ne sont pas fibr\'ees par des
hypersph\`eres.

\section{Introduction and statement of results}

In $3$-dimensional Euclidean space, the minimal surfaces known as
"Riemann minimal surfaces" belong to a one parameter family of
minimal surfaces which are embedded, have planar ends and are simply
periodic (i.e. are invariant under a discrete one parameter group of
translations). Moreover, in the quotient space, they have the
topology of a $2$-torus and have finite total curvature.

\medskip

These minimal surfaces have been discovered by Riemann in the $19$th
century and each element of this family is foliated by circles or
straight lines. In fact, up to some rigid motion and dilation, a
fundamental piece of any of these surfaces can be parameterized by
\[
X (t, \theta ) : = \left( a (t) + R  (t) \, \cos \theta, R  (t) \,
\sin \theta , t \right)
\]
for $(t, \theta) \in {\mathbb R} \times S^1$ in which case the
functions $a$ and $R$ are solutions of the following system of
first order nonlinear ordinary differential equations
\[
(\partial_t R)^2 +1 = \mu \, R^2 + R^4
\]
and
\[
\partial_t a  = R^2.
\]
where $\mu \in {\mathbb R}$ is a parameter. The surface
$\Sigma_{\mu}$ is invariant under some translation ${\bf d}_{\mu }
\in {\mathbb R}^3$, and, in the quotient space $\mathbb{R}^{3} / (
{\mathbb Z} \, {\bf d}_{\mu} )$, is topologically equivalent to a
torus $S^1 \times S^1$. In addition, in the quotient space, this
surface has finite total curvature.

\medskip

In this paper, we prove that this one parameter family of minimal
surfaces can be generalized to any dimension $n \geq 3$. More
precisely, we show that there exists a one parameter family of
minimal hypersurfaces in ${\mathbb R}^{n+1}$, which are embedded,
have infinitely many hyperplanar ends and are invariant under some
one parameter discrete group of translations.

\medskip

The canonical basis of $\mathbb{R}^{n+1}$ will be denoted by ${\bf
e}_{j}$, for $j=1, \ldots ,n+1$, and coordinates of $x\in {\mathbb
R}^n$ will be denoted by $(x^1, \ldots, x^n)$. In order to state our
result precisely, we introduce the subgroup ${\mathfrak G} \subset
O(n+1)$ which is generated by elements of the form
\[
R : = \begin{pmatrix}
  -1 & 0 & 0 \\
   0 & \bar R & 0 \\
   0 & 0 & -1 \\
\end{pmatrix}
\]
where $\bar R \in O(n-1)$.

\medskip

Our main result reads~:
\begin{theorem}
There exists a one parameter family of embedded minimal
hypersurfaces $(\Sigma_\e)_{\e \in (0, \e_0)}$ which have horizontal
hyperplanar ends and are simply periodic. These hypersurfaces are
invariant under the discrete group of translations ${\mathbb Z}\,
{\bf d}_\e$, where ${\bf d}_\e  = {\bf e}_1+ h_\e \, {\bf e}_{n+1}$
with $h_\e
>0$, and they are also invariant under the action of ${\mathfrak
G}$. In the quotient space ${\mathbb R}^{n+1} / ({\mathbb Z} \, {\bf
d}_\e )$ the hypersurface $\Sigma_\e$ is topologically equivalent to
$S^{n-1} \times S^1$ and has finite total curvature.
 \label{th:1}
\end{theorem}

When $n \geq 3$ and by opposition to the case of surfaces, it
follows from the result of W. C. Jagy \cite{Wil} that our
hypersurfaces are not foliated anymore by $(n-1)$-dimensional
spheres and hence it seems unlikely that these hypersurfaces could
be recovered by solving some system of nonlinear ordinary
differential equations. In fact, near any of its ends, the
hypersurface $\Sigma_\e$ we construct is close to the vertical graph
of the function
\[
x \in {\mathbb R}^n  \longrightarrow \frac{\e}{n-2} \, (
|x-x_*|^{2-n} - |x+x_*|^{2-n} )
\]
where $x_* =  (1, 0, \ldots, , 0)\in {\mathbb R}^n$ and one can
check that the level sets of this function are not spheres. As
already mentioned, near their ends, the minimal hypersurfaces we
construct are not exactly vertical graphs of these functions but
vertical graphs of some small perturbation these functions and it
turns out that, in dimension $n \geq 3$, the perturbations are small
enough so that the level sets are not spheres. This is in striking
contrast with the $2$-dimensional case where, the corresponding
$2$-dimensional construction  (leading to the construction of
Riemann minimal surfaces) yields a minimal surface which, close to
its ends, is close to the vertical graph of the function
\[
x \in {\mathbb R}^2  \longrightarrow \e \, ( \log \, |x+x_*| - \log
|x-x_*| )
\]
and this time one can check that the level sets of this function are
already circles.

\medskip

Let us emphasize that the hypersurfaces we construct do not describe
the full family of such hypersurfaces. Indeed, we only describe the
elements of this family when the translation period ${\bf d}_\e$ is
close to ${\bf e}_1$.

\medskip

In order to construct these hypersurfaces, the main observation is
that the moduli space of Riemann's minimal surfaces is one
dimensional (once the action of rigid motions and dilation has been
taken into account) and noncompact. In particular, one can
investigate the behavior of these surfaces close to one of the two
ends of the moduli space. It turns out that surfaces belonging to
one end of the moduli space, when they are properly rescaled, can be
understood as infinitely many parallel planes which are connected
together by small catenoidal necks. Hence, even though this would
not be worth the effort, these surfaces could be recovered using the
connected sum result by M. Traizet \cite{Tra}. This is this
connected sum result which allows us to describe part of the moduli
space of the $n$-dimensional analogues of Riemann's minimal
surfaces.

\medskip

This work compliments previous work which have been done to
generalize, in higher dimensions, some classical families of
minimal surfaces. For example, in \cite{Fak-Pac}, the minimal
$k$-noids, which are complete minimal surfaces with catenoidal
ends have been generalized to any dimension. This is also the case
for Scherk's second surfaces which have been generalized to any
dimension in \cite{Pac}.

\medskip

In section 2 we give the definition of the $n$-catenoid, which
generalizes the usual catenoid to any dimension. We then proceed
with a perturbation of the truncated $n$-catenoid to produce an
infinite dimensional family of minimal hypersurfaces which are
parameterized by their boundary data. Section 3 is devoted to the
perturbation of the hyperplane with two balls removed. Again we find
an infinite dimensional family of minimal hypersurfaces which are
parameterized by their boundary data on the boundaries of the two
excised balls. In Section 4, we explain how these infinite
dimensional families can be connected together to produce the
$n$-dimensional analogues of Riemann's minimal surfaces.

\section{The $n$-catenoid and minimal hypersurfaces close to it}

From now on, we assume that $n\geq 3$ is fixed. We recall some well
known facts concerning the definition and properties of the
$n$-catenoid $C$, a minimal hypersurface of revolution which
generalize in ${\mathbb R}^{n+1}$, the standard catenoid in
$3$-dimensional Euclidean space. We also give an rather explicit
expansion formula for the mean curvature of any hypersurface which
is close enough to $C$.

\medskip

The $n$-catenoid $C$ is a hypersurface of revolution about the
$x^{n+1}$-axis. It will be convenient to consider a parametrization
$\mathnormal{X}: \mathbb{R} \times
S^{n-1}\longrightarrow\,\mathbb{R}^{n+1}$ of $C$ for which the
induced metric  is conformal to the product metric on $\mathbb{R}
\times S^{n-1}$. This parametrization is given by
\begin{equation}
X(t,z) := \big( \varphi(t) \, z \, , \psi(t) \big),
\label{eq:catenoid C_0}
\end{equation}
where $t \in {\mathbb R}$, $z \in S^{n-1}$ and where the functions
$\varphi$ and $\psi$ are explicitly given by
\[
\varphi(t): = \left( \cosh ((n-1)t) \right)^\frac{1}{n-1} \qquad
\textrm{and}
 \qquad \psi(t) : = \displaystyle\int_{0}^{t}\varphi^{2-n} \, ds.
\]
It is easy to check that the induced metric on $C$ is given by
\[
g := \varphi^{2} \, (dt^2 +  g_{S^{n-1}})
\]
and, if the orientation of $C$ is chosen so that the unit normal
vector field is given by
\begin{equation}
{\bf n} : = \big( - \vph^{1-n} \, z , \partial_t \, \ln \vph\big),
\label{NormalGaussN0}
\end{equation}
then, the second fundamental form of $C$ is given by
\[
b  : ={\vph}^{2-n}\, \left ((1-n) \, dt^2  + g_{S^{n-1}} \right).
\]
From these expressions, it is easy to check that the hypersurface
parameterized by $X$ is indeed minimal.

\subsection{The Jacobi operator about the $n$-catenoid}

We now consider the hypersurfaces which can be parameterized as
normal graphs over $C$, namely they can be parameterized by
\begin{equation}
X_{w} := X + w \, {\bf n}
\label{eq:33}
\end{equation}
for some small (sufficiently smooth) function $w$. Let us denote by
$H(w)$ the mean curvature of the hypersurface parameterized by
$X_{w}$. The Jacobi operator, which is nothing but the linearized
mean curvature operator, appears in the second variation of the
$n$-volume functional. It is given by the general formula
\[
J : = \Delta_{g} + \arrowvert A \arrowvert^{2},
\]
where $\Delta_{g}$ denotes the Laplace-Beltrami operator and $A$
is the shape operator of the hypersurface. In the case of the
$n$-catenoid and in the above defined parametrization, the Jacobi
operator about $C$ is given by
\[
J = \varphi^{-n} \, \partial_t \, \left( \varphi^{n-2} \,
\partial_t \, \cdot \, \right) +  \varphi^{-2} \,
\Delta_{S^{n-1}} + n \, (n-1) \,  \varphi^{-2n}
\]
It turns out that it is easier (and equivalent) to study the mapping
properties of the conjugate operator $L$ which is defined by
\[
L : = \vph^{\frac{2+n}{2}} \, J \,  \vph^{\frac{2-n}{2}}.
\]
We have explicitly
\[
L : = \partial^2_t + \Delta_{{S^{n-1}}} - \left( \frac{_{n-2}}{^{2}}
\right)^{2} + \frac{_{n (3n-2)}}{^4} \, \vph^{2-2n}.
\]

The next Lemma is borrowed from \cite{Fak-Pac}. It explains the
structure of the expansion of the mean curvature operator $w
\longrightarrow H(w)$ in terms of the function $w$ and its
derivatives.
\begin{lemma} \cite{Fak-Pac}
The equation $H(w)=0$ is equivalent to
\begin{equation}
L \, w =  \vph^{\frac{2-n}{2}} Q _{2}
\big(\vph^{-\frac{n}{2}}w\big)\;+\;\vph^{\frac{n}{2}}Q _{3}
\big(\vph^{-\frac{n}{2}}w\big), \label{JacobiOperatorConj}
\end{equation}
where the operators $Q_{2}$ and $Q_3$ enjoy the following property :
There exists a constant $c>0$ such that for all $t\in\mathbb{R}$ and
for all $w_{1}, w_{2} \in
\mathcal{C}^{2,\alpha}([t-1,t+1]\times\Sn)$, we have
\begin{equation}
\|Q_{2}(w_{2})-Q_{2}(w_{1})\| _{\mathcal{C}^{0,\alpha}}\leq \,c \,
\big(\|w_{2}\|_{\mathcal{C}^{2,\alpha}}+\|w_{1}\|_{\mathcal{C}^{2,\alpha}}\big)
\,  \|w_{2}-w_{1}\|_{\mathcal{C}^{2,\alpha}}, \label{PropertiesQ2}
\end{equation}
and, provided $\|w_1\|_{\mathcal{C}^{2,\alpha}}+\|w_2
\|_{\mathcal{C}^{2,\alpha}} \leq 1$, we also have
\begin{equation}
\|Q_{3}(w_{2})-Q_{3}(w_{1})\| _{\mathcal{C}^{0,\alpha}}\leq \,c \,
\big(\|w_{2}\|_{\mathcal{C}^{2,\alpha}}+\|w_{1}\|_{\mathcal{C}^{2,\alpha}}\big)^{2}
\, \|w_{2}-w_{1}\|_{\mathcal{C}^{2,\alpha}}. \label{PropertiesQ3}
\end{equation}
Here all norms are understood on the domain of definition of the
functions. \label{lem:H(w)=0}
\end{lemma}
{\bf Proof~:} We recall the main lines of the proof of this critical
Lemma for the sake of completeness. We set
\[
\tilde{\bf n} : =  \varphi \, {\bf n} \qquad\mbox{and} \qquad
\tilde{w} : = \frac{w}{\varphi},
\]
so that the hypersurface parameterized by $X_w$ is also be
parameterized by $\tilde X_w : = X + \tilde{w} \, \tilde{\bf n}$.
Now, the first fundamental form $g_w$ of the hypersurface
parameterized by $\tilde X_w$ is explicitly given by
\[
\begin{array}{rlllll}
g_w & =  &  \displaystyle {\varphi}^2 \, (dt^2 + dz^i \, dz^j ) + 2
\, {\varphi}^{3-n}\, \tilde w \, ( (n-1) \, dt^2 +
dz^i \, dz^j)\\[3mm]
&  + &  \displaystyle 2 \, \varphi \, {\partial_t\varphi} \,  \tilde
w \, (
\partial_t \tilde w \, dt^2 +\partial_{z^i} \tilde w \, dt \, dz^i )
+ \varphi^{4-2n} \, \tilde w^2 \, ( n(n-2) \, dt^2 + dz^i \, dz^j )\\[3mm]
& + &  \displaystyle \varphi^2 \, ( (\tilde w^2 +(\partial_t
\tilde{w})^2 )\, dt^2 + 2 \, \partial_t \tilde{w} \,
\partial_{z^i} \tilde{w} \, dt \, dz^i + \partial_{z^i} \tilde w  \,
\partial_{z^j} \tilde w \,dz^i \, dz^j ),
\end{array}
\]
where $g_{S^{n-1}} = dz^i \, dz^j$ in local coordinates. Making use
of the expansion
\[
\mbox{det} \,  (I + B) = 1 + \mbox{Tr} \, B + \frac{_1}{^2} \left(
(\mbox{Tr} \, B)^2 - \mbox{Tr} (B^2)\right)+ {\cal O}(|B|^3)),
\]
and changing back $\tilde w $  into $w / \varphi$, we obtain the
expansion
\[
\sqrt {\mbox{det} \, g_w}  = \varphi^{n} + \frac{_1}{^2} \,
\varphi^{n-2} \, |\nabla w|^2 - \frac{_{n \, (n-1)}}{^2} \,
\varphi^{-n} \, w^2 + \varphi \, \tilde{Q}_3 \left( {\varphi}^{-1}
\, w \right) + \varphi^{n} \, \tilde{Q}_4\left( {\varphi}^{-1} \,
w\right)
\]
where $\tilde{Q}_3$ is homogeneous of degree $3$ and where
$\tilde{Q}_4$ collects all the higher order terms. The key point is
that the Taylor's coefficients of $\tilde{Q}_i$ are bounded
functions of $t$ and $z$ and so are the derivatives of any order of
these functions.

\medskip

The result then follows from the variational characterization of
minimal hypersurfaces as critical points of the functional
\[
{\cal E}( w )= \int \sqrt {\mbox{det} \, g_w}  \, ds \, dz
\]
It is easy to check that critical points of ${\cal E}$  are solution
of the nonlinear elliptic equation
\[
\partial_t(\varphi^{n-2} \partial_t w)+ \varphi^{n-2} \, \Delta_{S^{n-1}} \,w
+ n(n-1) \, \varphi^{-n} \, w = Q_2 \left(  {\varphi}^{-1} \,
w\right) + \phi^{n-1} \, Q_3 \left(  {\varphi}^{-1} \, w \right),
\]
where $Q_2$ is homogeneous of degree $2$ and where $Q_3$ collects
all the higher order terms. Again, the Taylor's coefficients of
$Q_i$ are bounded functions of $t$ and $z$ and so are the
derivatives of any order of these functions. To complete the proof,
it is enough to perform the conjugacy which was used to define $L$
starting from $J$. \hfill $\Box$

\medskip

Let us briefly comment of this result. The first estimate reflects
the fact that the operator $Q_2$ is a nonlinear second order
differential operator which is homogenous of degree $2$ in $w$ and
its derivatives, and
 has coefficients which are bounded functions of $t$. The second
estimate reflects the fact that the nonlinear operator $Q_{3}$ is
a nonlinear second order differential operator whose Taylor
expansion at $w=0$ does not involve any constant, linear nor
quadratic term and has coefficients which are bounded functions of
$t$.

\medskip

Observe that the $n$-catenoid is invariant under the action of the
group ${\mathfrak G}$ and if one looks for hypersurfaces which are
invariant under the action of ${\mathfrak G}$ then this amounts to
consider normal variations of the $n$-catenoid for some functions
$w$ which enjoy the following invariance property
\begin{equation}
w(- t, - z ) =  w ( t, z ) \qquad  \mbox{and} \qquad w( t, z ) = w
(t,  R \, z) \label{eq:sympr}
\end{equation}
for all $R \in O(n)$ of the form
\[
R : = \begin{pmatrix}
  1 & 0  \\
  0 & \bar R \\
\end{pmatrix}
\]
where $\bar R \in O(n-1)$. Clearly the Jacobi operator and its
conjugate preserve this invariance i.e. if a function $w$ satisfies
(\ref{eq:sympr}) then so does the function $L \, w $. Since the mean
curvature is invariant under the action of isometries, the nonlinear
operator which appears on the right hand side of
(\ref{JacobiOperatorConj}) also enjoys a similar invariance
property.

\subsection{Linear analysis about the n-catenoid}

We study the mapping properties of the conjugate Jacobi operator
$L$.

\medskip

Given $n \geq 2$, we denote by $\lambda_j = j \, (n-2+j)$, $j \in
{\mathbb N}$, the eigenvalues of the Laplace-Beltrami operator on
$S^{n-1}$ and we denote by $E_j$ the corresponding eigenspace. That
is
\[
\Delta_{S^{n- 1}} \phi = - \lambda_j \, \phi,
\]
for all $\phi \in E_j$.

\medskip

The {\em indicial roots} of $L$ describe the asymptotic behavior, at
infinity, of the solutions of the homogeneous problem
\[
L \, w = 0  \qquad \mbox{in} \qquad   \RS .  \\
\]
If $w \in\,\mathcal C^{2,\alpha}\big(\mathbb{R} \times S^{n-1}
\big)$ is solution of the homogeneous problem $\calL w = 0$ in
$\RS$, we may consider the eigenfunction decomposition of $w$ as
\[
w (t, z) = \sum_{j \in \mathbb N} w_j (t,z)
\]
where, for each $t\in {\mathbb R}$ the function $w_j(t, \cdot) \in
E_j$. Then the $E_j$-valued function $w_{j}$ is a solution of $L_{j}
\, v =0$ on $\mathbb{R}$ where the operator $L_{j}$ is defined by
\[
L_{j}:= \partial_{t}^{2}-  \left(\frac{_{n-2}}{^2} + j\right)^{2} +
\frac{_{n (3n-2)}}{^4} \, \vph^{2-2n}
\]
Since $\vph^{2-2n}$ tends exponentially to $0$ at $\pm \infty$, it
is easy to check that there are exactly two independent solutions of
the homogenous problem $L_{j} \, v =0$, which we denote by
$v_{j}^{+}$ and $v_{j}^{-}$ and which satisfy
\[
\lim_{t\longrightarrow + \infty} e^{-\gamma_{j} t} \,v_{j}^{+}(t) =1
\textrm{\quad\quad\quad  and  \quad\quad\quad}
\lim_{t\longrightarrow + \infty} e^{\,\gamma_{j} t} \, w_{j}^{-}(t)
=1
\]
where
\[
\delta_{j} : = \frac{_{n-2}}{^2}+ j,
\]
for all $j\,\in\,\mathbb{N} $. The real numbers $\pm \delta_{j}$ are
usually referred to as the {\em indicial roots} of $L$ at both
$+\infty$ and $-\infty$.

\medskip

We define the operator
\[
\Delta_{0} : = \partial_t^2+ \Delta_{S^{n-1}} - \left(
\frac{_{n-2}}{^2}\right)^2
\]
Observe that the indicial roots of $\Delta_0$ are  equal to the
indicial roots of the operator $L$.

\medskip

The solutions of the homogeneous problem $J \, w=0$ are usually
called  {\em Jacobi fields}. Some Jacobi fields which correspond to
explicit $1$-parameter family of minimal hypersurfaces to which $C$
belongs, are explicitly known. These Jacobi fields are obtained by
projecting over the normal vector field the Killing vector fields
associated to rigid motions and dilation. With slight abuse of
terminology we shall also refer to solutions of $L \, w = 0$ as
Jacobi fields. Then one has to multiply the Jacobi fields associated
to $J$ by $\varphi^{\frac{n-2}{2}}$ to obtain the expression of the
corresponding Jacobi fields for $L$, in doing so one takes into
account the fact that $L$ is conjugate to $L$. Using this receipt we
obtain the following Jacobi fields (for the operator $L$)~:

\begin{itemize}
\item[(i)] The function $\Phi^{0,-}  : = \varphi^{\frac{n-4}{2}} \, \partial_{t}\vph$,
which is associated to the translation of $C$ along its axis.

\item[(ii)] The function $\Phi^{0,+} : = \vph^{\frac{n-4}{2}} \, \big(\vph \,
\partial_{t}\psi - \psi \,
\partial_{t}\vph\big)$, which is associated to the dilation of $C$,

\item[(iii)] The function $\Phi^{1,-}_{\bf e} : = \varphi^{-\frac{n}{2}}
\, (z \, \cdot \, {\bf e})$ for ${\bf e} \in {\mathbb R}^n \times
\{0\}$, which is associated to the translation of $C$ along the
direction ${\bf e}$ orthogonal to its axis.

\item[(iv)] The function $\Phi^{1,+}_{\bf e} : = {\varphi}^{\frac{n-4}{2}} \,  (\psi \,
\partial_{t} \psi + \varphi \, \partial_{t} \varphi) \, (z \, \cdot  \, {\bf e})$ for ${\bf e}
\in {\mathbb R}^n \times \{0\}$,  which is associated to the
rotation of the axis of $C$ in a direction ${\bf e}$ orthogonal to
its axis.

\end{itemize}
These constitute $2 \, (n+1)$ linearly independent Jacobi fields.

\medskip

The operator $L$ does not satisfy the maximum principle. Indeed, one
checks that the Jacobi fields $\Phi^{1, -}_{\bf e}$ decay
exponentially at $\pm \infty$. Nevertheless, if the operator $L$ is
restricted to suitable subspace of functions, some version of the
maximum principle is still available. This is the content of the
following result whose proof can be found in \cite{Fak-Pac} (see
also \cite{Jle}).
\begin{proposition} \cite{Fak-Pac}
Assume that $ \delta < - \frac{n}{2}$. Let $w$ be a solution of
\[
L\,w =0 \textrm{\quad in \quad \,\,} {\mathbb R} \times S^{n-1},
\]
which is bounded by a constant times $(\cosh t)^{\delta}$. Then $w
\equiv 0$. \label{prop1maxim-princpl}
\end{proposition}
{\bf Proof~:} A simple proof of this result can be obtained as
follows. Proceed with the eigenfunction decomposition of $w$, the
solution of $L\, w =0$, so that $w  = \sum_j w_j$. Then $w_j$ is a
solution of $L_j \, w_j=0$ which is bounded by a constant times
$(\cosh t)^{\delta}$. When $j=0$ or $j=1$ then all solutions are
explicitly known and are described above. It is easy to check that
no solution is bounded  by a constant times $(\cosh t)^{\delta}$
unless it is identically equal to $0$ since we have chosen $\delta <
 - \frac{n}{2}$. Now, when $j \geq 2$ we write $w_j (t, z) = v_j (t)
 \, \phi_j(z)$ where $\phi_j \in E_j$.
As already observe the function $v_j$ being bounded by a constant
times $(\cosh t)^{\delta}$ for $\delta < - \frac{n}{2}$ has to decay
at infinity like $(\cosh t)^{-\delta_j}$. Then, we have
\[
\partial_{t}^{2} \, v_j -  \left(\frac{_{n-2}}{^2} + j\right)^{2}\, v_j +
\frac{_{n (3n-2)}}{^4} \, \vph^{2-2n} \, v_j =0
\]
and also (since $\Phi^{1, -}_{\bf e}$ are Jacobi fields)
\[
\partial_{t}^{2} \, \varphi^{-\frac{n}{2}} -  \left(\frac{_{n-2}}{^2} + 1\right)^{2}\, \varphi^{-\frac{n}{2}} +
\frac{_{n (3n-2)}}{^4} \, \vph^{2-2n} \, \varphi^{-\frac{n}{2}} =0
\]
For all $s \in {\mathbb R}$, we set $v_s : = \varphi^{-\frac{n}{2}}
- s \, v_j$. Using the above equations, we have
\begin{equation}
\partial_{t}^{2} \, v_s -  \left(\frac{_{n-2}}{^2} + j \right)^{2}\, v_s +
\frac{_{n (3n-2)}}{^4} \, \vph^{2-2n} \,  v_s  = - (j-1)
\,\left(\frac{_{n-2}}{^2} + j \right) \, \varphi^{-\frac{n}{2}}
\label{eq:cfacil} \end{equation} For all $s \in {\mathbb R}$, $v_s$
is positive near $\pm \infty$ (because the function $v_j$ tends to
$0$ at $\pm \infty$ much faster than the function
$\varphi^{-\frac{n}{2}}$). We choose $s$  to be the sup of the reals
for which $v_s \geq 0$. Then $v_s$ vanishes in ${\mathbb R}$ and at
this point, which is a minimum point for $v_s$, (\ref{eq:cfacil})
yields $\partial_{t}^{2} \, v_s  <0$. A contradiction. This
completes the proof of the result. \hfill $\Box$

\medskip

Given $k\,\in\,\mathbb{N} $, $\alpha\,\in(0,1)$ and $ \delta \, \in
\mathbb{R}$, the space $\mathcal{C}^{k,\alpha}_{\delta}({\mathbb R}
\times\Sn)$ is defined to be the space of functions $ w \, \in
\mathcal{C}^{k,\alpha} ({\mathbb R}\times \Sn) $ for which the
following norm
\[
\|w\|_{\mathcal{C}^{k,\alpha}_{\delta}({\mathbb R}\times\Sn)}:=\|
(\cosh t)^{-\delta} \, w\|_{\mathcal{C}^{k,\alpha}({\mathbb
R}\times\Sn)},
\]
is finite.

\medskip

We show in the next result that, provided the weight parameter
$\delta$ is suitably chosen and $\alpha \in (0,1)$ is fixed, one
can define a right inverse for $\mathcal L$.
\begin{proposition}
Assume that $\delta \in (\frac{n}{2} , \frac{n+2}{2})$ and $\alpha
\in (0,1)$ are fixed. Then, there exists a continuous operator
\[
G \, : \, \mathcal C^{0,\alpha}_{\delta} ({\mathbb R} \times \Sn)
\longrightarrow \mathcal C^{2,\alpha}_{\delta} ({\mathbb R}
\times\Sn)
\]
a right inverse for the operator $L$ such that, if the function $f$
satisfies (\ref{eq:sympr}) then so does $G \, (f)$. \label{prop:Surj
int}
\end{proposition}
\begin{dem}
The existence of $G$ follows from standard results and we refer to
\cite{Mel} and \cite{Maz} for a proof (see also \cite{Pac-Riv}).
According to Proposition~\ref{prop1maxim-princpl} and the
description of the geometric Jacobi fields, we see that the operator
\[
L  \, : \, \mathcal C^{2,\alpha}_{\delta} ({\mathbb R} \times \Sn)
\longrightarrow \mathcal C^{0,\alpha}_{\delta} ({\mathbb R}
\times\Sn)
\]
is injective for all $\delta < - \frac{n}{2}$. Hence according to
\cite{Mel} and \cite{Maz}, the operator $L$ is surjective for all
$\delta > \frac{n}{2}$ which is not an indicial root. This proves
the existence of a right inverse. Observe that, for $\delta >
\frac{n}{2}$ there is no uniqueness of the right inverse and in
order to define a right inverse which preserves (\ref{eq:sympr}) it
is enough to average over the orbit of the group and define
\[
G(f)(t,z)  = \frac{1}{2} \, \int_{O(n-1)} \left( \tilde G(f) (t,
(z^1,  \bar R \bar z)) + \tilde G(f) (-t, (-z^1, \bar R \bar
z))\right) \, d\sigma_{\bar R}
\]
where $z = (z^1, \bar z) \in S^{n-1}$ and where $\tilde G$ is {\em
any} right inverse for $L$. Here  $d\sigma_{\Bar R}$ is the standard
Haar measure on $O(n-1)$ (normalized so that the volume of $O(n-1)$
is equal to $1$).
\end{dem}

\medskip

The last result we will need is concerned with the Poisson
operator associated to the operator
$$
\Delta_{0} = \partial_{t}^2+\Delta_{\Sn}
-\big(\frac{_{n-2}}{^2}\big)^{2}
$$
which acts on functions defined on the cylinder ${\mathbb R}
\times S^{n-1}$. A similar result has already been proven in
\cite{Fak-Pac} but we give here a new very short self contained
proof.
\begin{lemma}
There exists a constant $ c =  c(n) > 0 $ such that for all
$h\,\in\mathcal C^{2,\alpha}(\Sn)$, which is
$L^2(S^{n-1})$-orthogonal to $E_0$ and $E_1$, there exists a unique
$w_h \in \mathcal C^{2,\alpha}([0,+\infty)\times\Sn)$ solution of
\[
\left\{\begin{array}{rllll} \Delta_{0} \, w_h &=&0&\textrm{\quad
in \quad} [0,+\infty)\times\Sn
\\[3mm]
w_h & = & h &\textrm{\quad on \quad}  \{0\}\times\Sn
\end{array}\right.
\]
which tends to $0$ as $t$ tends to $+\infty$. Furthermore,
\[
\| e^{\frac{n+2}{2}t} \, w_h \|_{\mathcal
C^{2,\alpha}([0,+\infty)\times\Sn)} \leq\, c \, \| h \|_{\mathcal
C^{2,\alpha}(\Sn)}.
\]
for some constant $c = c(n) >0$. \label{lem:HarmonicExten int}
\end{lemma}
\begin{dem}
We perform the eigenfunction decomposition of $h$
\[
h =\sum_{j \geq 2} h_{j}
\]
where $h_j \in E_j$. Then we have the explicit expression of $w_h$
given by
\[
w_h (t,z)=\sum_{j \geq 2}  e^{- \delta_j t} \, h_{j} (z),
\]
Using elliptic estimates together the fact that $h_j$ is an
eigenfunctions of $\Delta_{S^{n-1}}$, we get the rough estimate
\[
\| h_j \|_{L^\infty (S^{n-1})} \leq c \, (1+j)^{p_n} \, \| h_j
\|_{L^2 (S^{n-1})}
\]
for some exponent $p_n \geq 0$ which only depends on the dimension
$n$. Furthermore, using the fact that the dimension of $E_j$ grows
polynomially with $j$ \cite{Shu.book}, we get
\[
\|h_j \|_{L^2 (S^{n-1})} \leq c \, (1+j)^{q_n} \, \|h\|_{\mathcal
C^{2,\alpha}(\Sn)}.
\]
for some exponent $p_n \geq 0$ which only depends on the dimension
$n$. Collecting these, we conclude that
\[
e^{\frac{n+2}{2}t} \, |w_h (t,z)| \leq  c \, \left( \sum_{j \geq 2}
e^{(2- j) \, t} \, (1+j)^{p_n+q_n} \right)\, \|h\|_{\mathcal
C^{2,\alpha}(\Sn)}.
\]
It is easy to check that the series converges uniformly when $t \geq
1$. This provides the bound
\[
\sup_{t \geq 1 , z \in S^{n-1}} e^{\frac{n+2}{2}t} \, |w_h (t,z)|
\leq c \, \|h\|_{\mathcal C^{2,\alpha}(\Sn)}.
\]

Finally, applying the maximum principle we get \[\sup_{t \geq 0 , z
\in S^{n-1}} e^{\frac{n+2}{2}t} \, |w_h (t,z)| \leq c \,
\|h\|_{\mathcal C^{2,\alpha}(\Sn)}.\] The estimate for the
derivatives follows from Schauder's elliptic estimates.
\end{dem}

\medskip

Observe that if the function $h$ enjoys the following invariance
property
\begin{equation}
h( z ) =  h ( R \,z )\label{eq:symprh}
\end{equation}
for all $R \in O(n)$ of the form
\[
R : = \begin{pmatrix}
  1 &    0   \\
  0 & \bar R \\
\end{pmatrix}
\]
where $\bar R \in O(n-1)$, then so does the function $w_h$, namely
\begin{equation}
w_h (t, z ) =  w_h ( t, R \,z )\label{eq:symprhj}
\end{equation}
for all $R \in O(n)$ as above.

\subsection{Nonlinear analysis}

For all $\e \in (0,1)$, we define $t_{\e}>0$ by
\[
\varphi^{n-1}  (t_\e) = \e^{- \frac{{n}}{{3n-2}}}.
\]
and $r_\e >0$ is defined by
\[
r_\e : =  \e^{\frac{1}{n-1}} \, \varphi (t_\e) = \e^{\frac{2}{3n-2}}
\]
It will be convenient to slightly modify the normal vector field
${\bf n}$ on $C$ close into a transverse vector field ${\bf n}_{\e}$
which is defined by
\begin{equation}
{\bf n}_{\e} (t, z):= \chi_\e(t) \, {\bf n}(t,z) +  \mbox{sign} (t)
\, (1 - \chi_\e(t)) \, {\bf e}_{n+1} \label{NormalGaussNe}
\end{equation}
where $\chi_{\e}$ is a cutoff function equal to $0$ for $|t|>
t_{\e}-1$ and equal to $1$ when $|t| < t_{\e}-2$. Observe that, as
$\e$ tends to $0$, ${\bf n}_\e$ is a small perturbation of ${\bf
n}$. This is made quantitatively precise in the estimate
\[
|\nabla^{k}({\bf n}_{\e} \cdot {\bf n}-1)|\leq c_{k} \, \varphi^{2
-2 n}
\]
which holds for $|t|\geq t_\e-2$ and $z \in S^{n-1}$, for some
constant $c_k >0$ only depending on $k \in {\mathbb N}$.

\medskip

We now look for minimal hypersurfaces which are close to $C$ and
which can be parameterized by
\[
X _{\e,w} = X + \vph^{\frac{2-n}{2}} \, w \,  {\bf n}_\e,
\]
for some (small) function $w$. We also ask that these minimal
hypersurfaces are invariant under the action of ${\mathfrak G}$.

\medskip

Using the result of Lemma \eqref{lem:H(w)=0}, one can check that the
hypersurface parameterized by $X_{\e, w}$ is minimal if and only if
the function $w$ satisfies
\begin{equation}
L \, w = L_{\e} \, w + \vph^{\frac{2-n}{2}} \, Q_{2, \e} (
\vph^{-\frac{n}{2}}w) + \vph^{\frac{n}{2}} \, Q_{3,\e}
(\vph^{-\frac{n}{2}} \, w), \label{L(w)=Q_{e}(w)}
\end{equation}
This formula is not exactly identical to \eqref{JacobiOperatorConj}
since the normal vector field ${\bf n}$ has been modified into the
transverse vector field ${\bf n}_\e$. Observe that $L_\e \equiv 0$,
$Q_{2, \e} \equiv Q_2$ and $Q_{3, \e} \equiv Q_3$ when $|t|\leq
t_\e-2$ since ${\bf n} = {\bf n}_\e$ in this range. Moreover, it
follows from \eqref{NormalGaussNe} that the coefficients of the
linear second order operator $L_\e$ are bounded by a constant times
$\varphi^{2-2n}$ (in any ${\cal C}^{k, \alpha}$ topology). Indeed,
$\tilde L : = L -L_\e$ is the corresponding linearized mean
curvature operator when hypersurfaces close to $C$ are parameterized
as graphs over $C$ using the vector field ${\bf n}_\e$. It is easy
to check \cite{Fak-Pac} that, since $C$ has constant mean curvature
(equal to $0$) then
\[
\tilde L \, w  =  L \, ({\bf n}_\e \cdot {\bf n} \, w)
\]
The estimates of the coefficients of $L_\e$ follow at once from this
identity using (\ref{NormalGaussNe}). Finally, the operators
$Q_{2,\e}$ and $Q_{3,\e}$ enjoy properties which are similar to
those enjoyed by $Q_{2}$ and $Q_{3}$, uniformly for $\e \in (0,
1/2)$. Details can be found for example in \cite{Fak-Pac}.

\medskip

Given $k\,\in\,\mathbb{N} $, $\alpha\,\in(0,1)$ and $ \delta \,
\in \mathbb{R}$, the space
$\mathcal{C}^{k,\alpha}_{\delta}([-t_\e, t_\e] \times\Sn)$ is
defined to be the space of functions $ w \, \in
\mathcal{C}^{k,\alpha} ({\mathbb R}\times \Sn)$ which is endowed
with the norm :
\[
\|w\|_{\mathcal{C}^{k,\alpha}_{\delta}([-t_\e,
t_\e]\times\Sn)}:=\| (\cosh t)^{-\delta} \,
w\|_{\mathcal{C}^{k,\alpha}([-t_\e, t_\e] \times\Sn)}.
\]

Given $\tilde h \in{\cal C}^{2, \alpha}(S^{n-1})$ which is
$L^2(S^{n-1})$-orthogonal to $E_0$ and $E_1$, we define $w_{\tilde
h}$ to be the harmonic extension of $\tilde h$ in a half cylinder,
for the operator $\Delta_{0}$, which is given by
Lemma~\ref{lem:HarmonicExten int}. Then we set
\[ \tilde{w}_{\tilde
h} (t,z) := w_{\tilde h} (t_\e - t , z) + w_{\tilde h} (t_{\e} + t ,
-z )
\]
for all $(t,z) \in [-t_{\e},t_\e ]\times\Sn$. Granted the estimate
provided in Lemma~\ref{lem:HarmonicExten int}, it is easy to check
that there exists a constant $c = c(n)>0$ such that
\begin{equation}
\| \tilde w_{\tilde h}
\|_{\mathcal{C}^{2,\alpha}_{\frac{n+2}{2}}([-t_\e,
t_\e]\times\Sn)} \leq c \, \varphi^{-\frac{n+2}{2}}(t_\e) \, \|
\tilde h \|_{{\mathcal C}^{2, \alpha} (S^{n-1})} \label{eq:sdf}
\end{equation}
It will be convenient to define an extension operator
\[
{\mathcal E}_\e : \mathcal{C}^{0,\alpha}_{\delta}([-t_\e,
t_\e]\times\Sn) \longrightarrow
\mathcal{C}^{0,\alpha}_{\delta}({\mathbb R} \times\Sn)
\]
as follows~:

\begin{enumerate}

\item[(i)] For all $(t, z) \in (-t_\e, t_\e) \times S^{n-1}$, we set
\[
{\mathcal E}_\e (f) (t, z) = f (t, z).
\]
\item[(ii)] For all $(t, z) \in (-t_\e - 1, - t_\e) \times S^{n-1}$
we set
\[
{\mathcal E}_\e (f) (t, z) = \chi (t_\e -t) \, f ( -t_\e, z).
\]
\item[(iii)] For all $(t, z) \in (t_\e, t_\e +1 ) \times S^{n-1}$
we set
\[
{\mathcal E}_\e (f) (t, z) = \chi (t - t_\e) \, f (t_\e, z).
\]
\item[(iv)] And finally, for all $(t, z) \notin (-t_\e-1, t_\e+1) \times S^{n-1}$
we set
\[
{\mathcal E}_\e (f) (t, z) = f (t, z).
\]
\end{enumerate}
where $\chi : {\mathbb R} \longrightarrow [0,1]$ is a smooth cutoff
function identically equal to $0$ for $t \geq 1$ and identically
equal to $1$ for $t \geq 0$.

\medskip

Obviously there exists a constant $c = c(n, \delta) >0$ such that
\[
||| \, {\mathcal E}_\e \, ||| \leq c.
\]

We denote by $I_\e$ the canonical imbedding
\[
I_\e : \mathcal{C}^{2,\alpha}_{\delta}({\mathbb R} \times\Sn)
\longrightarrow \mathcal{C}^{2,\alpha}_{\delta}([-t_\e,
t_\e]\times\Sn)
\]

We fix $\delta\,\in (\frac{n}{2},\frac{n+2}{2})$ and we look for a
solution $w$ of \eqref{L(w)=Q_{e}(w)} which is defined in $[-t_\e,
t_\e]\times S^{n-1}$ and which can be decomposed as $w = \tilde
w_{\tilde h} + v$ where $v \in
\mathcal{C}^{2,\alpha}_{\delta}([-t_\e, t_\e]\times\Sn)$ is small.
Thanks to the result of Proposition~\ref{prop:Surj int}, we see that
it is enough to find $v \in {\cal C}^{2, \alpha}_\delta ([-t_\e,
t_\e]\times S^{n-1})$ solution of
\begin{equation}
v = A_{\e,\tilde h} (v)
,\label{eq:ptfix int}
\end{equation}
where we have defined
\[
A_{\e, \tilde h} (v) : = I_\e \circ G  \circ {\mathcal E}_\e  \,
\left( L_{\e} \, (\tilde w_{\tilde h} + v) + \vph^{\frac{2-n}{2}} \,
Q_{2, \e} ( \vph^{-\frac{n}{2}} ( \tilde w_{\tilde h} + v)) +
\vph^{\frac{n}{2}} \, Q_{3,\e} (\vph^{-\frac{n}{2}} \, ( \tilde
w_{\tilde h} + v)) - {\mathcal L} \, \tilde w_{\tilde h} \right)
\]
The existence of a solution $v$ to this fixed point problem will be
a consequence of the following Lemma and the application of a fixed
point theorem for contraction mapping.
\begin{lemma}
Given $\kappa >0$, there exists $\e_\kappa >0$ and $c_\kappa = c(n,
\delta, \kappa) >0$ such that, for all $\e \in (0, \e_\kappa)$, if
$\| \tilde h\|_{{\mathcal C}^{2, \alpha} (S^{n-1})} \leq \kappa \,
\e \, r_\e \, \varphi^{\frac{n}{2}} (t_\e)$. Then
\[
\| A _{\e, \tilde h} (0) \|_{\mathcal{C}^{2,
\alpha}_{\delta}([-t_{\e},t_{\e}]\times \Sn)} \, \leq c_\kappa \, \,
\e \, r_\e \, \varphi^{-1} (t_\e).
\]
and
\[
\| A _{\e, \tilde h}(v_2) - A _{\e, \tilde h}(v_1)
\|_{\mathcal{C}^{2,\alpha}_{\delta}([ -t_{\e},t_{\e} ] \times \Sn)}
\leq \frac{1}{2} \, \| v_2 - v_1
\|_{\mathcal{C}^{2,\alpha}_{\delta}([-t_{\e},t_{\e}]\times \Sn)}
\]
for all $v_1, v_2 \in {\mathcal
C}^{2,\alpha}_{\delta}([-t_{\e},t_{\e}] \times\Sn)$ satisfying
\[
\|v_i \|_{ {\mathcal C}^{2,\alpha}_{\delta}([-t_{\e},t_{\e}]
\times\Sn)} \leq 2 \, c_\kappa \, \varphi^{-1} (t_\e) \, \e \,
r_\e .
\]
\end{lemma}
\begin{dem}
We use (\ref{eq:sdf}) together with the properties of $L_\e$ to
get
\[
\| L_{\e} \, \tilde w_{\tilde h}  \|_{{\mathcal C}^{0,
\alpha}_\delta ([-t_\e, t_\e] \times \Sn) } \leq c \,
\varphi^{2-2n-\delta}(t_\e) \, \| \tilde h\|_{\mathcal
C^{2,\alpha}(\Sn)} \leq c_\kappa \, \e \, r_\e \,
\varphi^{\frac{4-3n}{2}-\delta } (t_\e).
\]
and also
\[
\| L_{\e} \, v  \|_{{\mathcal C}^{0, \alpha}_\delta ([-t_\e, t_\e]
\times \Sn) } \leq c \, \varphi^{2-2n}(t_\e) \, \|  v \|_{\mathcal
C^{2,\alpha}_{\delta} ([-t_\e, t_\e] \times \Sn)}.
\]

Next, observe that ${\mathcal L} - \Delta_0 =
\frac{_{n(3n-2)}}{^{4}} \, \varphi^{2-2n} $, hence
\[
\begin{array}{rllllll}
\| {\mathcal L} \, \tilde w_{\tilde h} \|_{{\mathcal C}^{0,
\alpha}_\delta ([-t_\e, t_\e] \times \Sn) } & = & \| ({\mathcal L}
-\Delta_0)\, \tilde w_{\tilde h} \|_{{\mathcal C}^{0, \alpha}_\delta
([-t_\e, t_\e] \times \Sn) } \\[3mm]
& \leq & c \, \vph^{-\frac{n+2}{2}}(t_{\e}) \, \| \tilde
h\|_{\mathcal C^{2,\alpha}(\Sn)} \\[3mm]
& \leq & c_\kappa \, \e \, r_\e \, \varphi^{-1} (t_\e). \end{array}
\]

Using the properties of $Q_{2, \e}$, we can estimate, for all $\e$
small enough (say $ \e \in (0, \e_\kappa)$),
\[
\| \vph^{\frac{2-n}{2}} \, Q_{2,\e} (  \vph^{-\frac{n}{2}} \, \tilde
w_{\tilde h} ) \|_{{\mathcal C}^{0, \alpha}_\delta ([-t_\e, t_\e]
\times \Sn) } \leq c \, \varphi^{-2-n} (t_\e) \, \|
 \tilde h\|_{\mathcal C^{2,\alpha}(\Sn)}^{2} \leq c_\kappa \, \e \,
 r_\e \, \varphi^{1-3n} (t_\e)
\]
and
\[
\begin{array}{rlllll}
\| \vph^{\frac{2-n}{2}} \, \left( Q_{2,\e} ( \vph^{-\frac{n}{2}} \,
(\tilde w_{\tilde h} +v_2) ) - Q_{2,\e} ( \vph^{-\frac{n}{2}} \,
(\tilde w_{\tilde h} + v_1) )\right) \|_{{\mathcal C}^{0,
\alpha}_\delta ([-t_\e, t_\e] \times \Sn) } \\[3mm]
\qquad \qquad \qquad \qquad \qquad \leq c \, \e \, r_\e \,
\varphi^{-1} (t_\e) \, \| v_2 - v_1 \|_{\mathcal
C^{2,\alpha}_\delta([-t_\e, t_\e] \times \Sn)}
\end{array}
\]
provided $v_1, v_2 \in {\mathcal C}^{2, \alpha}_\delta ([-t_\e,
t_\e] \times \Sn)$ satisfy the assumption of the statement.

\medskip

Similarly, using the properties of $\tilde{Q}_{3, \e}$, we can
estimate, for all $\e$ small enough (say $ \e \in (0, \e_\kappa)$,
\[
\|  \vph^{\frac{n}{2}} \, Q_{3, \e} ( \vph^{-\frac{n}{2}} \, \tilde
w_{\tilde h} )\|_{{\mathcal C}^{0, \alpha}_\delta ([-t_\e, t_\e]
\times \Sn) } \leq c \, \vph^{- n - \delta }(t_{\e}) \, \| \tilde
h\|_{\mathcal C^{2,\alpha}(\Sn)}^{3} \leq c_\kappa \, \e \, r_\e \,
\varphi^{\frac{12-9n}{2} -\delta} ,
\]
and
\[
\begin{array}{rllllll}
\| \vph^{\frac{n}{2}} \, \left( Q_{3,\e} ( \vph^{-\frac{n}{2}} \,
(\tilde w_{\tilde h} +v_2) ) - Q_{3,\e} ( \vph^{-\frac{n}{2}} \,
(\tilde w_{\tilde h} + v_1) )\right) \|_{{\mathcal C}^{0,
\alpha}_\delta ([-t_\e, t_\e] \times \Sn) } \\[3mm]
\qquad \qquad \qquad \qquad \qquad \leq c \, \e^2 \, r_\e^2 \,
\varphi^{n} (t_\e) \, \| v_2 - v_1 \|_{\mathcal C^{2,\alpha}_\delta
([-t_\e, t_\e] \times \Sn)} \end{array}
\] provided $v_1, v_2 \in {\mathcal C}^{2,
\alpha}_\delta ([-t_\e, t_\e] \times \Sn)$ satisfy the assumption of
the statement.

\medskip

The result follows at once from these estimates together with the
fact that the norms of $G$, $I_\e$ and ${\mathcal E}_\e$ are bounded
independently of $\e$.
\end{dem}

\medskip

Collecting the previous results, we conclude that, given $\kappa
>0$, there exists $\e_\kappa > 0$ such that, for all $\e \in (0,
\e_\kappa)$ the mapping $A_{\e, \tilde h}$ is a contraction from
\[
\{ v \in \mathcal C^{2,\alpha}_{\delta}([-t_{\e},t_{\e}] \times \Sn)
\quad : \quad \|v\|_{ \mathcal C^{2,\alpha}_{\delta}
([-t_{\e},t_{\e}] \times\Sn)}\leq 2 \, c_\kappa \, \e \, r_\e \,
\varphi^{-1} (t_\e) \big\}
\]
into itself and hence has a unique fixed point $v_{\e, \tilde h}$
in this set.

\medskip

The hypersurface parameterized by
\[
X _{\e, \tilde h} = \e^{\frac{1}{n-1}} \, \left( X  +  \,
\vph^{\frac{2-n}{2}} \, ( \tilde w_{\tilde h} + v_{\e, \tilde h})
\, {\bf n}_\e \right),
\]
for $(t, z) \, \in [- t_{\e}, t_{\e}] \times \Sn$ is a minimal
hypersurface which will be denoted by $C_{\e, \tilde h}$. This
produces an infinite dimensional family of minimal hypersurfaces
which are close to the piece of the catenoid $C$ which is the image
of $[-t_{\e}, t_{\e}] \times \Sn$ by $X$. This family is
parameterized by the boundary data $\tilde h$.

\medskip

Observe that, if one wants to produce a hypersurface which is
invariant under the action of ${\mathfrak G}$ it is enough to
restrict the subset of function $\tilde h \in \mathcal C^{2,\alpha}
(\Sn)$ which enjoy in addition the following invariance property
\begin{equation}
\tilde h( z ) =  \tilde h (R \, z) \label{eq:invpr}
\end{equation}
for all $R \in O(n)$ of the form
\[
R : = \begin{pmatrix}
  1 & 0  \\
  0 & \bar R \\
\end{pmatrix}
\]
for some $\bar R \in O(n-1)$.

\subsection{Local description of the hypersurface
$C_{\e, \tilde h}$ near its boundaries}

We first recall the asymptotic expansion of the parametrization of
the ends of the $n$-catenoid. Starting from the parametrization of
the $n$-catenoid which was given in \eqref{eq:catenoid C_0}, we
perform the change of variables
\[
x = \varphi (t) \, z \in {\mathbb R}^n
\]
The lower end of the $n$-catenoid can be parameterized as a vertical
graph over ${\mathbb R}^n \times \{0\}$ for some function $u$ which
can be expanded as
\[
u (x) = - d_0 + \frac{1}{n-2} \, |x|^{2-n} + {\cal O} (|x|^{4-3n}),
\]
for $|x|$ large enough, where
\[
d_0 : = \lim_{t\rightarrow +\infty} \psi ( t ).
\]
This expansion follows at once from the definition of $\varphi$ and
$\psi$. We refer to \cite{Fak-Pac} for the details.

\medskip

We now consider the $n$-catenoid which has been scaled by a factor
$\e^{\frac{1}{n-1}}$. Its lower end is can be parameterized by
\[
x \longrightarrow \e^{\frac{1}{n-1}} \, ( x ,  - d_0 +
\frac{_1}{^{n-2}} \, |x|^{2-n} + {\cal O} (|x|^{4-3n})),
\]
for $|x|$ large enough. Changing $\e^{\frac{1}{n-1}}$ into $x$, we
see that the lower end of the scaled $n$-catenoid can also be
parameterized by
\begin{equation}
x \longrightarrow (x ,  - \e^{\frac{1}{n-1}} \, d_0 + \frac{_\e
}{^{n-2}} \, |x|^{2-n} + {\cal O} (\e^3 \, |x|^{4-3n})),
\label{eq:lec}
\end{equation}

We apply the analysis of the previous section and collect the
results. Close to its lower boundary, the minimal hypersurface
$C_{\e , h}$ can be described as a vertical graph over an annulus in
the horizontal hyperplane $x^{n+1}=0$, this is the purpose of
changing the normal vector field ${\bf n}$ into ${\bf n}_\e$. To
make this precise, recall that we have defined
\[
r_\e : =  \e^{\frac{1}{n-1}} \, \varphi (t_\e) =
\e^{\frac{2}{3n-2}}
\]
and we write the hypersurface $C_{\e , h}$, close to its lower
boundary, as the graph over $\bar B(0, r_\e)- B(0, r_\e/2)$ for a
function $u_{\e, h}$.

\medskip

Given $h \in {\mathcal C}^{2, \alpha} (S^{n-1})$, we define
$\tilde h$ by
\[
\tilde h = \e^{-\frac{1}{n-1}} \, \varphi^{\frac{n-2}{2}} (t_\e)\, h
\]
We set
\[
u_{\e, h}^0(x) : =  - \e^{\frac{1}{n-1}} \, d_0 + \frac{\e}{n-2}
\, |x|^{2-n} + w_{h}^i (x/r_\e)
\]
where $w_{h}^i$ is the harmonic extension of $h$ in $B(0,1)$ and
we set
\begin{equation}
v_{\e, h} : = u_{\e, h} - u_{\e, h}^0 \label{eq:CauchyData1}
\end{equation}

Following the construction of the previous section, we obtain
the~:
\begin{lemma}
There exists $c =  c(n, \delta) > 0$ and, for all $\kappa >0$, there
exists $\e_\kappa >0$ such that, for all $h \in \mathcal
C^{2,\alpha}(\Sn)$ which is $L^2(S^{n-1})$-orthogonal to $E_0$ and
$E_1$ and which satisfies $\|h \|_{ \mathcal C^{2,\alpha}(\Sn)} \leq
\kappa\,\e\,r^{2}_{\e},$ we have
\[
\| v_{\e, h} (r_\e \cdot) \|_{\mathcal C^{2,\alpha} (\bar
B(0,1)-B(0, 1/2))} \leq c \, \e \, r^{2}_{\e}.
\]
In addition,
\[
\| (v_{\e, h_2} - v_{\e, h_1}) (r_\e \cdot) \|_{\mathcal
C^{2,\alpha} (\bar B(0,1)-B(0, 1/2))} \leq c \, \varphi^{\delta
-\frac{n+2}{2}} (t_\e)  \,  \| h_2 - h_1\|_{{\mathcal C}^{2, \alpha}
(S^{n-1})}.
\]
if $h_2, h_1 \in \mathcal C^{2,\alpha}(\Sn)$ are
$L^2(S^{n-1})$-orthogonal to $E_0$ and $E_1$ and which satisfy
\[
\|h_i \|_{ \mathcal C^{2,\alpha}(\Sn)} \leq \kappa \, \e \,
r^{2}_{\e}.
\]
\label{lem:CauchyData1}
\end{lemma}
The key and crucial point is that the constant $c >0$ does not
depend on $\kappa$. The proof of this estimate follows from a simple
but tedious computation following the steps of the construction of
$C_{\e,h}$. Observe that, when $h=0$, then the difference between
$u_{\e, 0}$ and $u_{\e, 0}^0$ comes from the term ${\cal O} (\e^3 \,
|x|^{4-3n})$, with $|x|=r_\e$, which appears in (\ref{eq:lec}) the
expansion of the lower end of the $n$-catenoid.

\medskip

When $h\neq 0$, there are many discrepancies to take into account.
The first comes from the fact that, by construction, the function
$\e^{\frac{1}{n-1}} \, \varphi^{\frac{2-n}{2}}\, \tilde w_{\tilde
h}$ is not exactly equal to the function $h$  when $t=t_\e$ but the
difference between these two functions is bounded by a constant
(depending on $\kappa$) times $\varphi^{-2-n} (t_\e) \, \e \,
r_\e^2$, and hence which is uniformly bounded by a constant
(independent of $\kappa$) times $ \e \, r_\e^2$ if $\e$ is taken
small enough. Next, one has to take into account the fact that, for
$t$ close to $t_\e$, the coordinates $(t, \theta)$ are not the usual
cylindrical coordinates $r = \e^{\frac{1}{n-1}} \, e^{s}$ in
${\mathbb R}^n$ and hence the normal graph of $\e^{\frac{1}{n-1}} \,
\varphi^{\frac{2-n}{2}} (t) \, \tilde w_{\tilde h}$ is not exactly
equal to the vertical graph of $w_h^i$. This induces in
\eqref{eq:CauchyData1} another discrepancy which is bounded by a
constant (depending on $\kappa$) times $\varphi^{2-2n} (t_\e) \, \e
\, r_\e^2$, and again is uniformly bounded by a constant
(independent of $\kappa$) times $ \e \, r_\e^2$ if $\e$ is taken
small enough. Finally, there is a term which comes from the
perturbation $v_{\e, \tilde h}$ solution of the nonlinear problem,
and this induces a discrepancy which is bounded by a constant
(depending on $\kappa$) times $\varphi^{\delta -\frac{n+2}{2}}
(t_\e) \, \e \, r_\e^2$, and since $\delta -\frac{n+2}{2} <0$, is
uniformly bounded by a constant (independent of $\kappa$) times $ \e
\, r_\e^2$ if $\e$ is taken small enough.

\section{Minimal hypersurfaces which are graphs over a hyperplane}

\subsection{The mean curvature for graphs}

Assume that a hypersurface is a vertical graph over the hyperplane
$x^{n+1}=0$ for some function $w$, i.e. this hypersurface is
parameterized by
\[
x \in {\mathbb R}^n \longrightarrow (x, w(x)) \in {\mathbb
R}^{n+1}.
\]
We recall that this hypersurface is minimal if and only if $w$ is a
solution of
\[
\mbox{div} \,\bigg(\frac{ \nabla \, w}{(1+| \nabla \,
w|^{2})^{\frac{1}{2}}}\bigg) =0
\]
It will be more convenient to write this equation as
\begin{equation}
 \Delta w = \frac{\nabla^{2} w \, (\nabla w,\nabla
w)}{1+|\nabla w|^{2}}. \label{eq:extPb}
\end{equation}

We will be interested in vertical graphs which are invariant under
the action of the group ${\mathfrak G}$ which has been defined in
the introduction. This amount to restrict our attention to functions
$w$ which enjoy the following invariance property
\begin{equation}
w( -x  ) =  - w ( x) \qquad  \mbox{and} \qquad w( x ) = w ( R \, x)
\label{eq:sympr2}
\end{equation}
for all $R \in O(n)$ of the form
\[
R : = \begin{pmatrix}
  1 &     0  \\
  0 & \bar R \\
\end{pmatrix}
\]
where $\bar R \in O(n-1)$. Again the Laplacian preserve this
invariance i.e. if a function $w$ satisfies (\ref{eq:sympr2}) then
so does the function $\Delta \, w $ and since the mean curvature is
invariant under the action of isometries, the nonlinear operator
which appears on the right hand side of (\ref{eq:extPb}) also enjoys
a similar invariance property.

\subsection{Linear analysis of the Laplacian in weighted spaces}

We set $x^* : = (1, 0, \ldots, 0) \in {\mathbb R}^n$ and we define
\[
{\mathbb R}^n_* : =  {\mathbb R}^n -\{x^*, -x^*\}
\]

Given $k\,\in\,\mathbb{N} $, $\alpha\,\in(0,1)$ and $\mu,
\nu\,\in\mathbb{R}$, the space $\mathcal C^{k,\alpha} _{\mu,\nu}
({\mathbb R}^n_*)$ to be the space of functions $w \in \mathcal
C^{k,\alpha}_{loc} ({\mathbb R}^n_*)$  for which the following norm
is finite
\[
\begin{array}{lllll}
\| w \|_{ \mathcal C^{k,\alpha}_{\mu, \nu} ({\mathbb R}^n_*)}  & : =
&  \sup_{s \in (0, 1/2)} s^{-\nu} \, \| w (s \, \cdot + x^*) \|_{
\mathcal C^{k,\alpha}(\bar B(0,2)- B(0,1))}  \\[3mm]
& + &  \sup_{s \in (0, 1/2)} s^{-\nu} \, \| w (s \, \cdot - x^*)
\|_{ \mathcal C^{k,\alpha}(B(0,2)- B(0,1))} \\[3mm]
& + & \| w \|_{ \mathcal
C^{k,\alpha}(\bar B(0, 4) - (B(x^*, 1/2)\cup B(-x^*, 1/2)))} \\[3mm]
& + &  \sup_{s \in (2, +\infty)} s^{-\mu} \, \| w (s \, \cdot ) \|_{
\mathcal C^{k,\alpha}(\bar B(0,2)- B(0,1))}.
\end{array}
\]
is finite. Therefore, the weight parameter $\nu$ controls the
behavior of the function $u$ near the points $\pm x^*$ and the
weight parameter $\mu$ controls its behavior at infinity.

\medskip

The following result follows from \cite{Mel} and \cite{Maz} but is
also a simple consequence of the maximum principle.
\begin{proposition}
Assume that $\mu ,  \nu \in (2-n,0)$ are fixed. Then, there exists a
constant $c =  c(n, \mu, \nu) >0$ and, for all $r \in (0, 1/2)$,
there exists a continuous operator
\[
\Gamma \, : \, \mathcal C^{0,\alpha}_{\mu-2,\nu-2} ({\mathbb R}^n_*)
\longrightarrow \mathcal C^{2,\alpha}_{\mu,\nu} ({\mathbb R}^n_*)
\]
such that, for all $f\,\in \mathcal
C^{0,\alpha}_{\mu-2,\nu-2}({\mathbb R}^n_*)$, the function $w:=
\Gamma (f)$ is a solution of
\[
\Delta \, w =  f
\]
in ${\mathbb R}^n_*$. In addition, if $f$ satisfies
(\ref{eq:sympr2}) then so does $\Gamma (f)$. \label{Prop:inv
Laplacien}
\end{proposition}
\begin{dem} As mentioned the existence of $\Gamma$ follows from
the results in \cite{Mel} and \cite{Maz} (see also \cite{Pac-Riv}).
Clearly the operator
\[
\Delta \, : \, \mathcal C^{2,\alpha}_{\mu ,\nu } ({\mathbb R}^n_*)
\longrightarrow \mathcal C^{2,\alpha}_{\mu-2,\nu-2} ({\mathbb
R}^n_*)
\]
is injective when $\nu > 2-n$ and $\mu <0$. Hence, according to the
results in \cite{Mel} and \cite{Maz}, the operator is surjective
when $\nu < 0$ and $\mu > 2-n$ are not indicial roots of the
Laplacian (namely are not of the form $j$ or $2-n-j$ for some $j \in
{\mathbb N}$. This completes the proof of the result.

\medskip

We also provide a simple proof of this result based on the maximum
principle. Observe that, for all $\lambda \in {\mathbb R}$, we have
\[
\Delta \, |x|^\lambda  =  \lambda \, (n-2+\lambda) \,
|x|^{\lambda-2}
\]
and that $\lambda \, (n-2+\lambda) <0$ precisely when $\lambda \in
(2-n, 0)$.

\medskip

Now, let us first assume that $f_1$ is supported in $B(0,4) - \{
x^*, -x^*\}$. The above observation implies that the function
\[
x \longrightarrow |x-x^*|^\nu + |x+x^*|^\nu
\]
can be used as a barrier function to prove both the existence $w_1$
solution of $\Delta w_1 =  f_1$ in ${\mathbb R}^n_*$ as well as the
pointwise bound
\[
|w_1 (x) |  \leq c  \,  \| f_1 \|_{{\mathcal C}^{0, \alpha}_{\mu-2,
\nu-2} ({\mathbb R}^n_*)} \, (|x-x^*|^\nu + |x+x^*|^\nu)
\]
for some constant $c = c(n,\nu) >0$.  However, since $f_1$ is
supported in $B(0,4) - \{x_*, -x_* \}$, $w_1$ is harmonic in
${\mathbb R}^n - B(0,4)$ and hence, it follows from the maximum
principle that the function $x \longrightarrow |x|^{2-n}$ can be
used as a barrier function to prove that
\[
|w_1 (x)| \leq \sup_{\partial B(0,4)} |w_1| \, 4^{n-2} \, |x|^{2-n}
\leq c \, \| f_1\|_{{\mathcal C}^{0, \alpha}_{\mu-2, \nu-2}
({\mathbb R}^n_*)}
\]
where $c = c(n, \nu) >0$.

\medskip

Finally, let us assume that $f_2$ is supported in ${\mathbb R}^n -
B(0, 2)$. Then, the above observation implies that the function
\[
x \longrightarrow |x- x^*|^\mu + |x + x^*|^\mu
\]
can be used as a barrier function to prove both the existence $w_2$
solution of $\Delta w_2 =  f_2$ in ${\mathbb R}^n_*$ as well as the
pointwise bound
\[
|w_2 (x) |  \leq c  \,  \| f_2 \|_{{\mathcal C}^{0, \alpha}_{\mu-2,
\nu-2} ({\mathbb R}^n_*)} \, (|x-x^*|^\mu + |x+x^*|^\mu)
\]
for some constant $c = c(n,\mu) >0$.  However, since $f_2$ is
supported in ${\mathbb R}^n - B(0,2)$, $w_2$ is harmonic in $B(0,2)$
and hence, it follows from the maximum principle that
\[
|w_2 (x)| \leq  \sup_{\partial B(0,2)} |w_2| \leq c \, \|
f_2\|_{{\mathcal C}^{0, \alpha}_{\mu-2, \nu-2} ({\mathbb R}^n_*)}
\]
where $c = c(n, \mu) >0$.

\medskip

The existence of $\Gamma (f)$ follows from these consideration,
first decomposing $f = f_1+f_2$  where $f_1$ is a function supported
in $B(0,4)$ and $f_2$ a function supported in ${\mathbb R}^n
-B(0,4)$. We obtain a function $\Gamma(f) = w_1 + w_2$ where $w_1$
and $w_2$ are defined as above. Collecting the above estimates, we
know that \[ |w (x)| \leq c \, \| f\|_{{\mathcal C}^{0,
\alpha}_{\mu-2, \nu-2} ({\mathbb R}^n_*)} \, ( |x- x^*|^\nu + |x +
x^*|^\nu)\] in $B(0, 4)$ and
\[
|w (x)|\leq c \, \| f\|_{{\mathcal C}^{0, \alpha}_{\mu-2, \nu-2}
({\mathbb R}^n_*)} \, |x|^\mu \] in ${\mathbb R}^n- B(0, 4)$, for
some constant $c = c(n, \mu, \nu) >0$. Once the existence and the
pointwise control of $w$ have been obtained, the estimates for the
derivatives of $w$ follow from Schauder's estimates. This completes
the proof of the result.
\end{dem}

\medskip

The following Lemma is the counterpart of
Lemma~\ref{lem:HarmonicExten int}. Observe that, this time we do
not impose any constraint on the boundary data and the proof is
again a simple application of the maximum principle.
\begin{lemma} Assume that
$\alpha \in (0,1)$ is fixed. For all $h \in
\mathcal{C}^{2,\alpha}(S^{n-1})$, there exists $\bar w_h \in
\mathcal C^{2,\alpha}_{2-n} (\mathbb{R}^n - B(0,1))$ satisfying :
\[
\left\{\begin{array}{rllll} \Delta \bar w_h & = & 0 &\textrm{\quad
in \quad}\mathbb{R}^n -
\bar B(0,1)\\[3mm]
       \bar w_h & = & h &\textrm{\quad on\quad}\partial B(0,1)
\end{array} \right.
\]
Furthermore,
\[
\|\bar w_h \|_{{\mathcal C}^{2,\alpha}_{2-n}(\mathbb{R}^n - B(0,1))}
\leq c \, \|h \|_{\mathcal{C}^{2,\alpha}(S^{n-1})}.
\]
for some constant $c = c(n) >0$. \label{lem:HarmonicExtent ext}
\end{lemma}
\begin{dem}
The proof of this result follows at once from the observation that
$x \longrightarrow |x|^{2-n}$ can be used as a barrier function to
prove both the existence of $\bar w_h$ and the pointwise bound
\[
|\bar w_h (x) |\leq \sup_{\partial B(0,1)} |h| \, |x|^{2-n}.
\]
We then apply Schauder's estimates to get the relevant estimates for
the derivatives of $\bar w_h$.
\end{dem}

\medskip

Observe that if the function $h$ enjoys the invariance property
(\ref{eq:symprh}), namely
\[
h( z ) =  h ( R \,z )
\]
for all $R \in O(n)$ of the form
\[
R : = \begin{pmatrix}
  1 &    0   \\
  0 & \bar R \\
\end{pmatrix}
\]
where $\bar R \in O(n-1)$, then so does the function $\bar w_h$,
namely
\begin{equation}
\bar w_h ( x ) =  \bar w_h ( R \, x )\label{eq:symprhj}
\end{equation}
for all $R \in O(n)$ as above.

\subsection{Nonlinear analysis}

Recall that we have defined
\[
r_\e : = \e^{\frac{2}{3n-2}}.
\]
For all $r \in (0,1)$, we also define
\[
D_r : = \mathbb{R}^{n} - ( \bar B(x^* , r) \cup \bar B(-x^*, r)).
\]
Given $k\,\in\,\mathbb{N} $, $\alpha\,\in(0,1)$ and $\mu,
\nu\,\in\mathbb{R}$, we define the space $\mathcal C^{k,\alpha}
_{\mu,\nu} (D_r)$ to be the space of functions $w \in \mathcal
C^{k,\alpha}_{loc} (D_r)$  for which the following norm
\[
\begin{array}{lllll}
\| w \|_{ \mathcal C^{k,\alpha}_{\mu, \nu} (D_r)}  & : = &
 \sup_{s \in (r, 1/2)} s^{-\nu} \, \| w (s \, \cdot + x^*) \|_{
\mathcal C^{k,\alpha}(\bar B(0,2)- B(0,1))}  \\[3mm]
& + &  \sup_{s \in (r, 1/2)} s^{-\nu} \, \| w (s \, \cdot - x^*)
\|_{\mathcal C^{k,\alpha}(\bar B(0,2)- B(0,1))} \\[3mm]
& + & \| w \|_{ \mathcal C^{k,\alpha}(\bar B(0, 4) - (B(x^*, 1/2)\cup B(-x^*, 1/2)))} \\[3mm]
& + &  \sup_{s \in (2, +\infty)} s^{-\mu} \, \| w (s \, \cdot ) \|_{
\mathcal C^{k,\alpha}(\bar B(0,2)- B(0,1))}.
\end{array}
\]
is finite.

\medskip

The remaining of the analysis parallels what we have already done
to perturb the truncated rescaled $n-$catenoid. Given $\bar h \in
{\cal C}^{2, \alpha} (S^{n-1})$, we use the result of
Lemma~\ref{lem:HarmonicExtent ext} to define the function
\[
\hat w_{\bar h}  : =  \bar w_{\bar h} ( (\cdot - x^*)/ r_\e )  -
\bar w_{\bar h} ( - (\cdot + x^*)/ r_\e ).
\]
Using the estimate provided by Lemma~\ref{lem:HarmonicExtent ext},
we conclude that
\[
\| \hat w_{\bar h} \|_{{\mathcal C}^{2, \alpha}_{2-n,2-n}(D_{r_\e})}
\leq  c \, r_\e^{n-2} \, \| \bar h\|_{{\cal C}^{2,
\alpha}(S^{n-1})}.
\]
for some constant $c = c(n) >0$.

\begin{remark}
Observe that the decay at infinity of $\hat w_{\bar h}$ can be
slightly improved and we have
\[
\| \hat w_{\bar h} \|_{{\mathcal C}^{2, \alpha}_{1-n,2-n}(D_{r_\e})}
\leq  c \, r_\e^{n-2} \, \| \bar h\|_{{\cal C}^{2,
\alpha}(S^{n-1})}.
\]
This follows from the fact that the function $\bar w_{\bar h}$ can
be decomposed as
\[
\bar w_{\bar h} (x) = a_{\bar h} \, |x|^{2-n} + \bar w_{\bar h}' (x)
\]
where
\[
|a_{\bar h} |+ \| \bar w_{\bar h}' \|_{{\mathcal C}^{2,
\alpha}_{1-n} ({\mathbb R}^n - B(0,1))} \leq c \, \|\bar h
\|_{{\mathcal C}^{2, \alpha} (S^{n-1})}.
\]
\end{remark}

We define, for all $\rho \in {\mathbb R}$ and $\e
>0$ the function
\[
w_{\e, \rho} (x) : = (\rho - 2^{1-n}\, \e) \, x^1 + \frac{\e}{n-2}
\, ( |x-x^*|^{2-n}- |x+x^*|^{2-n}).
\]
We look for a solution of \eqref{eq:extPb} of the form
\[
w := w_{\e , \rho} + \hat w_{\bar h} + v,
\]
where $v$ is a small function.

\medskip

It will be convenient to define an extension operator
\[
\bar {\mathcal E}_\e : \mathcal{C}^{0,\alpha}_{\mu, \nu}(D_{r_\e})
\longrightarrow \mathcal{C}^{0,\alpha}_{\mu, \nu}({\mathbb R}^n_*)
\]
as follows~:

\begin{enumerate}

\item[(i)] For all $x \in D_{r_\e}$, we set
\[
\bar {\mathcal E}_\e (f) (x) = f (x).
\]
\item[(ii)] For all $ B(x^*, 2r_\e)- \bar B(x^*, r_\e)$
we set
\[
\bar {\mathcal E}_\e (f) (x) = \bar \chi (|x-x^*|/ r_\e) \, f ( x^*
+ r_\e \, (x-x_*)/|x-x^*|).
\]
\item[(iii)] For all $ B(-x^*, 2r_\e)- \bar B(-x^*, r_\e)$
we set
\[
\bar {\mathcal E}_\e (f) (x) = \bar \chi (|x+x^*|/ r_\e) \, f ( -
x^* + r_\e \, (x+x_*)/|x+x^*|).
\]
\item[(iv)] And finally, for all $x \in B(x^*, r_\e) \cup B(-x^*, r_\e)$
we set
\[
\bar {\mathcal E}_\e (f) (0) = 0.
\]
\end{enumerate}
where $\bar \chi : {\mathbb R} \longrightarrow [0,1]$ is a smooth
cutoff function identically equal to $1$ for $s \geq 2$ and
identically equal to $0$ for $s \leq 1$.

\medskip

Obviously there exists a constant $c = c(n , \mu , \nu) >0$ such
that
\[
||| \, \bar {\mathcal E}_\e \, ||| \leq c.
\]

We denote by $\bar I_\e$ the canonical imbedding
\[
\bar I_\e : \mathcal{C}^{2,\alpha}_{\delta}({\mathbb R}^n_*)
\longrightarrow \mathcal{C}^{2,\alpha}_{\delta}(D_{r_\e})
\]

\medskip

We fix $\mu, \nu \in (2-n, 0)$. Using the result of
Proposition~\ref{Prop:inv Laplacien}, we rephrase this problem as a
fixed point problem. It is now enough to find $v \in \mathcal
C^{k,\alpha}_{\mu,\nu}(D_{r_\e})$ solution of
\begin{equation}
v  = B_{\e,\rho , \bar h} (v)
\label{eq:ptfix ext}
\end{equation}
where
\[
B_{\e,\rho , \bar h} (v) : = \bar I_\e \circ \Gamma \circ \bar
{\mathcal E}_\e \,  (\Xi (w_{\e, \rho} + \hat w_{\bar h} + v)) .
\]
and where we have set
\[
\Xi (w) : = \frac{\nabla^{2}w(\nabla w,\nabla w)}{1+|\nabla w|^{2}}
\]

The existence of a fixed point for $B_{\e,\rho , \bar h}$ relies on
the following~:
\begin{lemma}
There exists $c =  c(n, \mu, \nu ) >0$ and for all $\kappa>0$ there
exists $\e_{\kappa}>0$ such that for all $\varepsilon\,\in
(0,\varepsilon_{\kappa})$, for all $\rho \in {\mathbb R}$ and for
all $\bar h \in {\mathcal C}^{2, \alpha} (S^{n-1})$ satisfying
\[
r_\e \, |\rho| + \| \bar h\|_{{\mathcal C}^{2, \alpha} (S^{n-1})}
\leq \kappa \, \e \, r_\e^2 ,
\]
we have
\[
\| B_{\e,\rho , \bar h} (0) \|_{{\mathcal C}^{2, \alpha}_{\mu, \nu}
(D_{r_\e})} \leq c \, \e \, r_\e^{2-\nu}.
\]
Moreover,
\[
\|B_{\e,\rho, \bar h}(v_1) - B_{\e,\rho, \bar h }(v_2)\|_{ \mathcal
C^{2,\alpha}_{\mu,\nu}(D_{r_\e})} \leq \frac{1}{2} \, \| v_1 - v_2
\|_{\mathcal C^{2,\alpha}_{\mu, \nu}(D_{r_\e})}
\]
for all $v_1, v_2  \in {\mathcal C^{2,\alpha}_{\mu,\nu}(D_{r_\e})}$
satisfying
\[
\|v_i\|_{\mathcal C^{2,\alpha}_{\mu,\nu}(D_{r_\e})} \leq 2 \, c \,
\e \, r_\e^{2-\nu}.
\]
\label{lem:1}
\end{lemma}
\begin{dem}
The first estimate follows from the result of
Proposition~\ref{Prop:inv Laplacien} together with the estimate
\[
\| \Xi (w_{\e, \rho} + \hat w_{\bar h} )\|_{\mathcal
C^{0,\alpha}_{\mu-2,\nu-2}(D_{r_\e})} \leq c \, \e \, r_\e^{2-\nu}
\]
which follows from the construction of $w_{\e, \rho} + \hat w_{\bar
h}$. The second estimate follows from
\[
\| \Xi (w_{\e, \rho} + \hat w_{\bar h} + v_2) - \Xi (w_{\e, \rho} +
\hat w_{\bar h} +v_1 )\|_{{\mathcal
C}^{0,\alpha}_{\mu-2,\nu-2}(D_{r_\e})} \leq c \, r_\e^n \, \| v_2
-v_1\|_{{\mathcal C}^{2,\alpha}_{\mu,\nu}(D_{r_\e})}
\]

Details are left to the reader.
\end{dem}

\medskip

Collecting the previous results, we conclude that, given $\kappa
>0$, there exists $\e_\kappa > 0$ such that, for all $\e \in (0,
\e_\kappa)$ the mapping $K_{\e, \rho , \bar h}$ is a contraction
from
\[
\{ v \in \mathcal C^{2,\alpha}_{\mu, \nu }(D_{r_\e}) \quad : \quad
\|v\|_{ \mathcal C^{2,\alpha}_{\mu, \nu} (D_{r_\e}) } \leq 2 \, c \,
\e \, r_\e^2 \big\}
\]
into itself and hence has a unique fixed point $w_{\e, \rho , \bar
h}$ in this set. We define
\[
\bar u_{\e, \rho , \bar h} : = w_{\e , \rho}  + \hat w_{\bar h} +
w_{\e, \rho, \bar h}.
\]
The hypersurface parameterized by $ x \longrightarrow (x, \bar
u_{\e, \rho , \bar h}(x))$ for $x \in D_{r_\e}$ is a minimal
hypersurface and will be denoted by $\Sigma_{\e, \rho, \bar h}$.

\medskip

The important fact is that the constant $c>0$ which appears in
Lemma~\ref{lem:1} does not depend on $\kappa$ provided $\e$ is
chosen small enough.

\medskip

If one looks for minimal hypersurfaces which are invariant under the
action of the group ${\mathfrak G}$, it is enough to retrict our
attention to the set of functions $\bar h \in {\mathcal C}^{2,
\alpha} (S^{n-1})$ which are invariant under
\[
\bar h(z) =   \bar h (R z)
\]
for all $R \in O(n)$ of the form
\[
R : = \begin{pmatrix}
  1 &    0   \\
  0 & \bar R \\
\end{pmatrix}
\]
where $\bar R \in O(n-1)$.

\subsubsection{Local description of the hypersurface $\Sigma_{\e, \rho, \bar h}$ near its
boundaries}

We would like to analyze $\bar u_{\e, \rho , \bar h}$ close to
$x^*$. To this aim, we define $y : =  x-x^*$ and the function
\[
\bar u^0_{\e , \rho , \bar h} (y) : = \rho -  \frac{n}{n-2} \,
2^{1-n} \, \e + \frac{\e}{n-2} \, |y|^{2-n} + \rho \, y^1 + \bar
w^e_{\bar h} (y/r_\e)
\]
where
\[
w^e_{\bar h} : =  \bar w_{\bar h}
\]
is the harmonic extension defined in Lemma 3.1. We also define
\[
\bar v_{\e, \rho, \bar h} : = \bar u_{\e , \rho , \bar h} - \bar
u^0_{\e , \rho , \bar h}
\]
Following the construction of $\bar u_{\e , \rho , \bar h}$, we
obtain~:
\begin{lemma}
There exists $c =  c(n, \mu, \nu)>0$  and for all $\kappa>0$ there
exists $\e_\kappa >0$  such that, for all $\e \in (0,\e_\kappa )$,
for all $\rho \in {\mathbb R}$ and for all $\bar h  \in  \mathcal
C^{2,\alpha}(\Sn)$ satisfying
\[
r_\e \, |\rho| + \|\bar h \|_{ \mathcal C^{2,\alpha}(\Sn)}\leq
\kappa \, \e \, r^{2}_{\e}
\]
we have
\[
\|\bar v_{\e, \rho, \bar h}  (x^* + r_\e \,  \cdot ) \|_{ \mathcal
C^{2,\alpha}(_bar B(0,2)-B(0,1))} \leq c \, \e \, r^{2}_{\e}.
\]
Moreover,
\[
\|(\bar v_{\e, \rho_2, \bar h_2} - \bar v_{\e, \rho_1, \bar h_1})
(x^* + r_\e \, \cdot ) \|_{ \mathcal C^{2,\alpha}(\bar
B(0,2)-B(0,1))} \leq c \, r_\e^{n} \, (r_\e \, |\rho_2 - \rho_1|+ \|
\bar h_2 - \bar h_1 \|_{{\mathcal C}^{2, \alpha} (S^{n-1})})
\]
for all
\[
r_\e \, |\rho_i| + \|\bar h_i \|_{ \mathcal C^{2,\alpha}(\Sn)}\leq
\kappa \, \e \, r^{2}_{\e}.
\]
 \label{lem:CauchyData2}
\end{lemma}
Again, the important fact is that the constant $c > 0$ in the first
estimate does not depend on $\kappa$.

\section{The connected sum construction}

We fix $\kappa$ large enough and apply the results of the previous
sections.

\medskip

Assume that we are given $h, \bar h \in {\mathcal C}^{2,
\alpha}(S^{n-1})$ satisfying
\[
\| h\|_{{\mathcal C}^{2, \alpha} (S^{n-1})} \leq \kappa \, \e \,
r_\e^2 \qquad \mbox{and} \qquad \| \bar h \|_{{\mathcal C}^{2,
\alpha} (S^{n-1})} \leq \kappa \, \e \, r_\e^2.
\]
We decompose
\[
h = h_0 +h_1 + h^\perp
\]
where $h_0 \in E_0$, $h_1 \in E_1$ and  $h^\perp$ is
$L^2(S^{n-1})$-orthogonal to $E_0$ and $E_1$. We further assume that
both $h$ and $\bar h$ satisfy \eqref{eq:invpr}. In particular, this
implies that $h_1 \in \mbox{Span} \{  x\cdot {\bf e}_1 \}$ and if we
decompose
\[
\bar h = \bar h_0 + \bar h_1 + \bar h^\perp
\]
this also implies that $\bar h_1 \in \mbox{Span} \{  x\cdot {\bf
e}_1 \}$.

\medskip

Granted the above decomposition, we choose from now on
\[
\rho : = - r_\e^{-1} \, h_1
\]
and consider the hypersurface $\Sigma_{\e, \rho , \bar h }$ which
has been constructed in the previous section. Next, we define
\[
t : = \left( \e^{\frac{1}{n-1}} \, d_0 + \rho - \frac{n}{n-2} \,
2^{1-n} \, \e + h_0  \right)
\]
and consider ${\bf e}_1 + t \, {\bf e}_{n+1}  + C_{\e, h^\perp}$ the
hypersurface defined in section  \S 2.4 and which has been
translated by ${\bf e}_1 + t \, {\bf e}_{n+1}$. The lower boundary
of ${\bf e}_1 + t \, {\bf e}_{n+1}  + C_{\e, h^\perp}$ and the
"upper boundary" of $\Sigma_{\e, \rho , \bar h }$ are close on to
the other and we will now show that, for all $\e$ small enough, it
is possible to find $\bar h$ and $h$ in such a way that the union of
$\Sigma_{\e, \rho , \bar h }$ and ${\bf e}_1 + t \, {\bf e}_{n+1} +
C_{\e, h^\perp}$ is a ${\mathcal C}^1$ hypersurface. Since this
hypersurface has piecewise mean curvature equal to $0$ and is
${\mathcal C}^1$, regularity theory then implies that it is a smooth
minimal hypersurface with two boundaries and one end asymptotic to
$x\longrightarrow (x, \rho \, x^1)$.

\medskip

To complete the construction of the generalized Riemann's minimal
hypersurface, it will remain to first apply a suitable rotation so
that the end of the hypersurface $\Sigma_{\e, \rho , \bar h } \cup
({\bf e}_1 + t \, {\bf e}_{n+1} + C_{\e, h^\perp})$ which is
asymptotic to $x\longrightarrow (x, \rho \,  x^1)$ becomes
horizontal and next to extend this hypersurface so that it becomes a
sinply periodic hypersurface (which is invariant under the action of
${\mathfrak G}$).

\medskip

In order to produce a $\mathcal C^{1}$ hypersurface, we consider the
two summands as vertical graphs over annular regions in the
hyperplane $x^{n+1}=0$ and ask that the Cauchy data of these two
graphs coincide. This condition can be translated into the following
set of equations
\begin{equation}
\left\{\begin{array}{rlllll} u_{\e, h^\perp}(r_\e \cdot) + t & = &
\bar u_{\e, \rho, \bar h} (x^* + r_\e \cdot) \\[3mm]
\partial_r u_{\e, h^\perp}(r_\e \cdot)  & = & \partial_r \bar
u_{\e, \rho, \bar h} (x^* + r_\e \cdot)
\end{array}\right.
\end{equation}
on $S^{n-1}$. However, given the expansions of $u_{\e, h^\perp}$
and $\bar u_{\e, \rho , \bar h}$ this is equivalent to solve
\begin{equation}
\left\{\begin{array}{lllllll} w^i_h - w^e_{\bar h}   & = &  \bar
v_{\e, \rho, \bar h} (x^*+ r_\e
\cdot)-  v_{\e, h^\perp} (r_\e  \cdot)  \\[3mm]
\partial_{r} (w^i_h - w^e_{\bar h}) & = &  \partial_r (\bar v_{\e,
\rho, \bar h} (x^*+ r_\e \cdot)-  v_{\e, h^\perp} (r_\e  \cdot))
\end{array}\right.\label{eq:Cond C^1}
\end{equation}
on $S^{n-1}$.

\medskip

To proceed with, we recall the following result \cite{Maz-Pac}
\begin{lemma}
The mapping
\[
\begin{array}{rclclll}
\mathcal P :& \mathcal C^{2,\alpha}(\Sn) &
\longrightarrow&\mathcal
C^{1,\alpha}(\Sn) \\[3mm]
& h  &\longmapsto    & \partial_{r} \, (w_{h}^e - w_{h}^i)
\end{array}
\]
is an isomorphism.
\end{lemma}

Using this result, the solvability of \eqref{eq:Cond C^1} reduces
to a fixed point problem which can be written as
\[
(h, \bar h) = S_\e (h, \bar h)
\]
It follows from the estimates of Lemma~\ref{lem:CauchyData1} and
Lemma~\ref{lem:CauchyData2} that
\[
\| S_\e (h, \bar h)\|_{({\cal C}^{2, \alpha}(S^{n-1}))^2} \leq c_0
\, \e \, r_\e^2
\]
for some constant $c_0 >0$ which does not depend on $\kappa$
provided $\e$ is small enough. In addition
\[
\| S_\e (h_2, \bar h_2) -S_\e (h_1, \bar h_1) \|_{({\cal C}^{2,
\alpha}(S^{n-1}))^2} \leq \frac{1}{2} \, \| (h_2-h_1, \bar h_2 -
\bar h_1)\|_{({\cal C}^{2, \alpha}(S^{n-1}))^2}
\]
provided $\e$ is chosen small enough.

\medskip

To conclude, we choose $\kappa = 2 \, c_0$ and use a fixed point
Theorem for contraction mappings which will ensure the existence of
at least one fixed point for the mapping $S_\e$ in
\[
\left\{ (h, \bar h) \in ({\mathcal C}^{2, \alpha} (S^{n-1}))^2
\quad : \quad \| (h , \bar h)\|_{({\mathcal C}^{2, \alpha}
(S^{n-1}))^2} \leq \kappa \, \e\, r_\e^2 \right\}
\]
provided $\e$ is chosen small enough, say $\e \in (0, \e_0]$. This
completes our proof of the existence of a fixed point for $S_\e$
and hence the existence of the Riemann minimal hypersurface in any
dimension.


\medskip


\begin{thebibliography}{99}

        \bibitem{Fak-Pac} S. Fakhi and F. Pacard, {\em Existence of
complete minimal hypersurfaces with finite total curvature,}
        Manuscripta Math. {\bf 103}, (2000), 465-512.

        \bibitem{Gil-Tru.book} D. Gilbarg and N.S. Trudinger, {\em
Elliptic Partial Differential Equations of second order}
        Springer, (2001).

        \bibitem{Jle} M. Jleli, {\em Constant mean curvature
hypersurfaces}, PhD Thesis, University of Paris
        12 (2004).

\bibitem{Maz} R. Mazzeo, {\em Elliptic theory of edge operators I.}
Comm. in PDE No 16, 10 (1991) 1616-1664.

        \bibitem{Maz-Pac} R. Mazzeo and F. Pacard {\em Constant scalar
curvature metrics with isolated
        singulaties}, Duke Math. J. {\bf 99}, (1999), 353-418.

        \bibitem{Meek-Per-Ros} W. Meeks, J. Perez and A. Ros, {\em
Uniqqueness of the Riemann minimal examples}
        Invent. Math. {\bf 131}, (1998), 107-132.

        \bibitem{Mel} R. Melrose, {\em The Atiyah-Patodi-Singer index
theorem}. Research Notes in Mathematics, xiv, 377 p. (1993).

\bibitem{Pac} F. Pacard, \emph{Higher dimensional Scherk's
hypersurfaces} J. Math. Pures Appl., IX. Ser. 81, No.3, (2002)
241-258.

\bibitem{Pac-Riv} F. Pacard and T. Rivi\`ere, {\em Linear and nonlinear
aspects of vortices : the Ginzburg Landau model} Progress in
Nonlinear Differential Equations, 39, Birk\"auser (2000).


\bibitem{Shu.book} M.A. Shubin, {\em Pseudodifferential
operators and spectral theory}, Springer, (1987).

\bibitem{Tra} M. Traizet, \emph{Adding handles to Riemann
minimal examples} J. Inst. Math. Jussieu 1, No.1, (2002) 145-174 .

\bibitem{Wil} W. C. Jagy {\em Minimal hypersurfaces foliated by
spheres}, Michigan Math. J. {\bf 38}, (1991), 255-270.

  \end{thebibliography}
\end{document}